# UNIQUENESS FOR DIFFUSIONS DEGENERATING AT THE BOUNDARY OF A SMOOTH BOUNDED SET

By Dante DeBlassie

*Texas A&M University*

For continuous $\gamma, g : [0,1] \to (0, \infty)$, consider the degenerate stochastic differential equation
$$dX_t = [1 - |X_t|^2]^{1/2} \gamma(|X_t|) \, dB_t - g(|X_t|) X_t \, dt$$
in the closed unit ball of $\mathbb{R}^n$. We introduce a new idea to show pathwise uniqueness holds when $\gamma$ and $g$ are Lipschitz and $\frac{g(1)}{\gamma^2(1)} > \sqrt{2} - 1$. When specialized to a case studied by Swart [*Stochastic Process. Appl.* **98** (2002) 131–149] with $\gamma = \sqrt{2}$ and $g \equiv c$, this gives an improvement of his result. Our method applies to more general contexts as well. Let $D$ be a bounded open set with $C^3$ boundary and suppose $h : \overline{D} \to \mathbb{R}$ Lipschitz on $\overline{D}$, as well as $C^2$ on a neighborhood of $\partial D$ with Lipschitz second partials there. Also assume $h > 0$ on $D$, $h = 0$ on $\partial D$ and $|\nabla h| > 0$ on $\partial D$. An example of such a function is $h(x) = d(x, \partial D)$. We give conditions which ensure pathwise uniqueness holds for
$$dX_t = h(X_t)^{1/2} \sigma(X_t) \, dB_t + b(X_t) \, dt$$
in $\overline{D}$.

**1. Introduction.** For a long time much has been known about uniqueness for one-dimensional stochastic differential equations (SDEs) with singular coefficients. The diffusion coefficient can be non-Lipschitz and degenerate; the drift can be singular and involve local time. See the survey (in Section 4) of Engelbert and Schmidt (1991), as well as the references there. In contrast, the higher-dimensional situation is understood less. Recent work in this direction includes the articles of Athreya, Barlow, Bass and Perkins (2002), Bass and Perkins (2002) and Swart (2001, 2002).

Athreya, Barlow, Bass and Perkins (2002) and Bass and Perkins (2002) study weak uniqueness for
$$dX_t^i = \sum_k \sqrt{2 X_t^i} \sigma_{ik}(X_t) \, dB_t^k + b^i(X_t) \, dt, \qquad i = 1, \ldots, n,$$

in the positive orthant in $\mathbb{R}^n$, where $b$ and $\sigma$ satisfy suitable nonnegativity and regularity conditions. This problem is interesting because the diffusion matrix is degenerate and non-Lipschitz and the boundary of the state space is not smooth.

Swart (2002) addressed both weak uniqueness and pathwise uniqueness for the SDE

$$(1.1) \qquad dX_t = \sqrt{2(1-|X_t|^2)}\, dB_t - cX_t\, dt$$

in the closed unit ball $E$ in $\mathbb{R}^n$. As above, the diffusion matrix is degenerate and non-Lipschitz. He proved weak uniqueness holds when $c \geq 0$ and $n \geq 1$. Standard methods yield pathwise uniqueness in dimension $n=1$ for all $c \geq 0$ and also in dimensions $n \geq 2$, provided $c=0$ or $c \geq 2$. The case $0 < c < 2$ for $n \geq 2$ is much trickier. Swart used a clever method to prove pathwise uniqueness for $c \geq 1$. Rotational invariance of (1.1) played a large role in the argument. Also, with the explicit form of the coefficients, Swart was able to exploit the resulting explicit form of the local time on the boundary. In this article we study a slightly more general form of (1.1) in the closed unit ball $E$ of $\mathbb{R}^n$:

$$(1.2) \qquad dX_t = [1-|X_t|^2]^{1/2} \gamma(|X_t|)\, dB_t - g(|X_t|) X_t\, dt,$$

where $\gamma, g : [0,1] \to (0, \infty)$.

We introduce a new technique yielding a theorem, which specialized to (1.1), improves Swart's result.

THEOREM 1.1. *Suppose $\gamma, g : [0,1] \to (0, \infty)$ are Lipschitz continuous with $\frac{g(1)}{\gamma(1)^2} > \sqrt{2} - 1$. Then pathwise uniqueness holds for (1.2).*

REMARK 1.1. In the context of (1.1), we have $\gamma(1) = \sqrt{2}$ and $g(1) = c$. Hence, the condition $\frac{g(1)}{\gamma(1)^2} > \sqrt{2} - 1$ becomes $c > 2(\sqrt{2} - 1) \approx 0.828$. This improves Swart's condition $c \geq 1$.

REMARK 1.2. Since the process $1-|X_t|^2$ is an autonomous one-dimensional diffusion, a change of space and time can be used to prove existence of a solution to (1.2). The idea is much like that used in the proof of Theorem 1.2.

It is natural to ask if the power $\frac{1}{2}$ in (1.2) can be changed to $r > 0$. When $r \geq 1$, the coefficients are Lipschitz and it is well known that pathwise uniqueness holds. When $r \in (\frac{1}{2}, 1)$, if the process starts within the open unit ball, then the boundary is unattainable [see the last chapter in Breiman (1968)] and, again, standard results yield pathwise uniqueness. If the process starts on the boundary, our method can be used to show pathwise uniqueness holds in this case too; more on this at the end of Section 3. Finally, when $r \in (0, \frac{1}{2})$, our method does not seem to work and we do not know if pathwise uniqueness holds. To see that *pathwise uniqueness* is the issue, in Section 4 we outline the proof of the following theorem. The technique is standard.

THEOREM 1.2. *If $\gamma, g : [0, 1] \to (0, \infty)$ are continuous, then for any $r \in (0, \frac{1}{2})$, weak uniqueness holds for*

$$dX_t = [1 - |X_t|^2]^r \gamma(|X_t|) \, dB_t - g(|X_t|) X_t \, dt,$$
$$X_0 = x \in E$$

*in the closed unit ball $E$ of $\mathbb{R}^n$.*

Now we explain the idea behind the proof of Theorem 1.1. For solutions $X$ and $\widetilde{X}$ to (1.2) with the same Brownian motion, the usual idea for proving pathwise uniqueness is to compute $d|X - \widetilde{X}|^2$, show the integrands of the resulting terms involving $dt$ are bounded by $|X - \widetilde{X}|^2$, then appeal to Gronwall's inequality. But due to the non-Lipschitz nature of the diffusion coefficient in (1.2), $d|X - \widetilde{X}|^2$ has a $dt$ term whose integrand $I$ is positive and singular in the sense that $\frac{I}{|X - \widetilde{X}|^2}$ is unbounded. This precludes the use of Gronwall's inequality. Swart's idea is to look at

$$W = (Y^{1/2} - \widetilde{Y}^{1/2})^2 + |X - \widetilde{X}|^2,$$

where $Y = 1 - |X|^2$ and $\widetilde{Y} = 1 - |\widetilde{X}|^2$. Here $d(Y^{1/2} - \widetilde{Y}^{1/2})^2$ gives rise to a negative singular term which, under the condition $c \geq 1$, compensates for the positive singular term in $d|X - \widetilde{X}|^2$. Our idea is to use

$$W = (Y^p - \widetilde{Y}^p)^2 + |X - \widetilde{X}|^2$$

for suitable $p \in (\frac{1}{2}, 1)$. For this choice of $p$ there will be an extra positive singular term in $d(Y^p - \widetilde{Y}^p)^2$ not occurring in Swart's work. The critical observation is under the condition $\frac{g(1)}{\gamma^2(1)} > \sqrt{2} - 1$, this new positive singular term can also be absorbed into the negative singular term. This is a bit surprising because creating more positive singular terms does not seem to be a good idea initially.

To simplify the exposition, we have chosen to concentrate on (1.2) in the closed unit ball. But our technique applies to more general situations, since we do not rely on explicit properties of the local time on the boundary to prove Theorem 1.1. Indeed, we now state a more general version of the theorem.

Let $D \subseteq \mathbb{R}^n$ ($n \geq 2$) be a bounded open set such that for some $\varphi \in C^3(\mathbb{R}^n)$,

$$D = \{x \in \mathbb{R}^n : \varphi(x) > 0\},$$
$$\partial D = \{x \in \mathbb{R}^n : \varphi(x) = 0\},$$
$$|\nabla \varphi| > 0 \text{ on } \partial D.$$

Suppose $h: \overline{D} \to \mathbb{R}$ satisfies

(1.3)
$$h > 0 \text{ on } D,$$
$$h = 0, \ |\nabla h| > 0 \text{ on } \partial D,$$
$$h \text{ is Lipschitz on } D,$$
$$h \text{ is } C^2 \text{ with Lipschitz second derivatives, all on a neighborhood of } \partial D.$$

An example of such a function is $h(x) = d(x, \partial D)$.

Consider the SDE

$$(1.4) \qquad dX_t = [h(X_t)]^{1/2} \sigma(X_t) \, dB_t + b(X_t) \, dt$$

in the closed set $\overline{D}$, where $B_t$ is a Brownian motion in $\mathbb{R}^n$, $\sigma = (\sigma_{ij})$ is an $n \times n$ matrix and $b$ is an $n$-dimensional vector, both Lipschitz on $\overline{D}$. Assume

$$(1.5) \qquad a(x) = \sigma(x)\sigma^*(x)$$

is strictly positive definite for $x \in \overline{D}$:

$$(1.6) \qquad \langle a(x)\xi, \xi \rangle > 0, \qquad x \in \overline{D}, \xi \in \mathbb{R}^n \setminus \{0\},$$

where $\langle \cdot, \cdot \rangle$ is the usual Euclidean inner product. We also assume there is a neighborhood $N$ of $\partial D$ such that

$$(1.7) \qquad b = g \frac{\nabla h}{|\nabla h|} + \beta \qquad \text{on } N,$$

where

(1.8)
$$g > 0 \text{ and Lipschitz on } N,$$
$$\beta \text{ is Lipschitz on } N,$$
$$\langle \beta, \nabla h \rangle = 0 \text{ on } N.$$

Then $g$ is uniformly bounded below away from 0. We say $f(x)$ is a Lipschitz function of $h$ if for some constant $C > 0$,

$$|f(x) - f(y)| \le C|h(x) - h(y)|.$$

Equivalently, $f = \bar{f} \circ h$ for some Lipschitz $\bar{f}$.

THEOREM 1.3. *In addition to conditions* (1.3) *and* (1.5)–(1.8), *suppose* $g|\nabla h|$ *and* $\langle a\nabla h, \nabla h \rangle$ *are Lipschitz functions of $h$ on a neighborhood of $\partial D$. Then automatically* $\alpha = \frac{2g\nabla h}{\langle a\nabla h, \nabla h \rangle}|_{\partial D}$ *is constant. If* $\alpha > \sqrt{2} - 1$, *then pathwise uniqueness holds for* (1.4).

The method of proof is like that for Theorem 1.1. Please note the condition requiring $g|\nabla h|$ and $\langle a\nabla h, \nabla h\rangle$ to be Lipschitz functions of $h$ is rather restrictive. For instance, the hypotheses do not cover a simple nonrotationally symmetric equation proposed by Swart:
$$dX_t = \sqrt{2(1-|X_t|^2)}\,dB_t + c(\theta - X_t)\,dt$$
in $E$, where $\theta \in \mathbb{R}^n$ is constant.

The article is organized as follows. In Section 2 we compute $d[1-|X|^2]^p$ for $p \in (\frac{1}{2}, 1)$ and show $X$ spends zero Lebesgue time on the boundary; the latter is needed for the proof of Theorem 1.1. In Section 3 we prove Theorem 1.1 and discuss the proof of Theorem 1.3, as well as the case $r > \frac{1}{2}$ mentioned after Remark 1.2. The proof of Theorem 1.2 is outlined in Section 4; the proof consists of standard methods. In Section 5 we present some open questions. The last section consists of the proof of a technical result used in the proof of Theorem 1.1.

**2. The differential of powers of $1 - |X_t|^2$.** Let $X$ be any solution to (1.2), where $g$ and $\gamma$ are continuous. For any $p > 1 - \frac{g(1)}{\gamma^2(1)}$, by continuity, choose $\varepsilon(p) > 0$ such that

(2.1) $$p > 1 - \frac{g(u)}{\gamma^2(u)}, \qquad u \in (1-\varepsilon(p), 1].$$

For any process $R$ and $\delta > 0$, define
$$\tau_\delta(R) = \inf\{t \geq 0 : R_t = \delta\}.$$

NOTATION. In the sequel we will write
$$dR = a\,dB + b\,dt,$$
for $t \leq \eta$, to mean
$$R_{t\wedge\eta} = R_0 + \int_0^{t\wedge\eta} a(s)\,dB(s) + \int_0^{t\wedge\eta} b(s)\,ds.$$

Here is our result on powers of $1-|X|^2$. We suppress the explicit dependence of $\gamma$ and $g$ on $|X|$. Notice no boundary local time terms appear—this is why we require $p > g(1)/\gamma^2(1)$.

LEMMA 2.1. *Let $p \in (\frac{1}{2}, 1)$ satisfy $p > \frac{g(1)}{\gamma^2(1)}$ and suppose $\varepsilon = \varepsilon(p)$ is from (2.1). Then for $\tau = \tau_\varepsilon(1-|X|^2)$ and $|X_0|^2 > 1 - \varepsilon$, the process $Y = 1 - |X|^2$ satisfies*
$$dY^p = -2pY^{p-1/2}\gamma \sum_j X_j\,dB_j$$
$$+ 2p|X|^2 Y^{p-1} I_{Y>0}[g + (p-1)\gamma^2]\,dt - np\gamma^2 Y^p\,dt$$

*for $t \leq \tau$.*

PROOF. For $m \geq 1$, choose $\ell_m \in C(\mathbb{R})$ with $\ell_m \geq 0$, $\mathrm{supp}\,\ell_m \subseteq [\frac{1}{m+1}, \frac{1}{m}]$ and $\int \ell_m(t)\,dt = 1$. Then

$$k_m(t) = \int_0^t \int_0^s \ell_m(u)\,du\,ds$$

satisfies $k_m \in C^2(\mathbb{R})$,

$$k_m(t) \to t \vee 0 \quad \text{as } m \to \infty, \text{ uniformly on } \mathbb{R},$$
$$k_m \equiv 0 \quad \text{on a neighborhood of } 0,$$
$$0 \leq k'_m \leq 1,$$
$$k'_m \to I_{(0,\infty)} \quad \text{as } m \to \infty,$$
$$k''_m \geq 0.$$

Now for $y = 1 - |x|^2$,

$$\frac{\partial}{\partial x_i} k_m(y^p) = -2px_i y^{p-1} k'_m(y^p),$$

$$\frac{\partial^2}{\partial x_i^2} k_m(y^p) = -2p y^{p-1} k'_m(y^p) + 4p(p-1) x_i^2 y^{p-2} k'_m(y^p)$$
$$+ 4p^2 x_i^2 y^{2p-2} k''_m(y^p).$$

Hence, for $t \leq \tau$,

$$d[k_m(Y^p)] = -2p Y^{p-1/2} \gamma k'_m(Y^p) \sum_j X_j\,dB_j$$
$$+ 2p|X|^2 Y^{p-1} k'_m(Y^p)[g + (p-1)\gamma^2]\,dt$$
(2.2)
$$- pn\gamma^2 Y^p k'_m(Y^p)\,dt$$
$$+ 2p^2 |X|^2 Y^{2p-1} \gamma^2 k''_m(Y^p)\,dt.$$

We are going to show the integrated forms of the first three terms on the right-hand side converge to their analogs with $k'_m(Y^p)$ replaced by $I_{Y>0}$ and the integrated form of the last term converges to 0. To this end, for any $t > 0$,

$$E\left[\left|\int_0^{t\wedge\tau} 2p Y^{p-1/2} \gamma [k'_m(Y^p) - I_{Y>0}] \sum_j X_j\,dB_j\right|^2\right]$$

(2.3)
$$= 4p^2 E\left[\int_0^{t\wedge\tau} Y^{2p-1} \gamma^2 [k'_m(Y^p) - I_{Y>0}]^2 |X|^2\,ds\right]$$

$$\leq 4p^2 \sup(\gamma^2) \varepsilon^{2p-1} E\left[\int_0^{t\wedge\tau} [k'_m(Y^p) - I_{Y>0}]^2\,ds\right]$$

(since $t \leq \tau$ and $|X_0|^2 > 1 - \varepsilon$ imply $Y_t \leq \varepsilon$). Properties of $k_m$ and dominated convergence show

(2.4) $\qquad$ left-hand side (2.3) $\to 0 \qquad$ as $m \to \infty$.

Also,

$$E\left[\left|\int_0^{t\wedge\tau} pn\gamma^2 Y^p[k'_m(Y^p) - I_{Y>0}]\, ds\right|\right]$$

(2.5) $$\leq pn\varepsilon^p \sup(\gamma^2) E\left[\int_0^{t\wedge\tau} |k'_m(Y^p) - I_{Y>0}|\, ds\right]$$

$$\to 0 \qquad \text{as } m \to \infty,$$

by dominated convergence. Finally,

(2.6) $\qquad E[|k_m(Y^p_{t\wedge\tau}) - Y^p_{t\wedge\tau}|] \to 0 \qquad$ as $m \to \infty$

by uniform convergence of $k_m(t)$ to $t \vee 0$.

Looking at the integrated form of (2.2) and using (2.4)–(2.6), we see

$$\int_0^{t\wedge\tau} 2p|X|^2 Y^{p-1} k'_m(Y^p)[g + (p-1)\gamma^2]\, ds + \int_0^{t\wedge\tau} 2p^2|X|^2 Y^{2p-1}\gamma^2 k''_m(Y^p)\, ds$$

must converge in $L^1$ as $m \to \infty$. Clearly, the second integrand is nonnegative and by (2.1), the first integrand is too. Hence, Fatou's lemma yields

$$E\left[\int_0^{t\wedge\tau} 2p|X|^2 Y^{p-1} I_{Y>0}[g + (p-1)\gamma^2]\, ds\right] < \infty,$$

and then, by dominated convergence,

(2.7) $$E\left[\left|\int_0^{t\wedge\tau} 2p|X|^2 Y^{p-1}[k'_m(Y^p) - I_{Y>0}][g + (p-1)\gamma^2]\, ds\right|\right]$$
$$\to 0 \qquad \text{as } m \to \infty.$$

Now we can let $m \to \infty$ in the integrated form of (2.2) and use (2.4)–(2.7) to end up with

(2.8) $$\begin{aligned} dY^p = {} & -2pY^{p-1/2}\gamma \sum_j X_j\, dB_j \\ & + 2p|X|^2 Y^{p-1} I_{Y>0}[g + (p-1)\gamma^2]\, dt \\ & - pn\gamma^2 Y^p\, dt \\ & + d\varphi_t^{(p)}, \end{aligned}$$

for $t \leq \tau$, where $\varphi_t^{(p)}$ is continuous and nondecreasing in $t$. The conclusion of the lemma will follow once we prove $\varphi_t^{(p)} \equiv 0$.

First we show $\varphi_t^{(p)}$ can change only when $Y_t = 0$:

$$\text{(2.9)} \qquad \int_0^{t\wedge\tau} I_{Y>0}\,d\varphi^{(p)} = 0.$$

By (2.8) and Itô's formula [since $k_m(0) = 0$],

$$dk_m(Y^p) = -k'_m(Y^p)2pY^{p-1/2}\gamma\sum_j X_j\,dB_j$$
$$+ k'_m(Y^p)2p|X|^2Y^{p-1}[g + (p-1)\gamma^2]\,dt$$
$$- k'_m(Y^p)pn\gamma^2 Y^p\,dt$$
$$+ k'_m(Y^p)\,d\varphi_t^{(p)}$$
$$+ k''_m(Y^p)2p^2 Y^{2p-1}\gamma^2|X|^2\,dt.$$

Compare with (2.2) to see we must have

$$\int_0^{t\wedge\tau} k'_m(Y^p)\,d\varphi^{(p)} = 0.$$

Let $m \to \infty$ and use dominated convergence to get

$$\int_0^{t\wedge\tau} I_{Y>0}\,d\varphi^{(p)} = 0,$$

as claimed.

To finish, choose $q > \frac{1}{2}$ such that $q < p$ and

$$q > 1 - \frac{g(u)}{\gamma^2(u)}, \qquad u \in (1 - \varepsilon(p), 1].$$

Then the derivation leading to (2.8) holds with $p$ replaced by $q$ and the analogue of (2.8) is valid. By an extension of Itô's formula to $C'$ functions [Rogers and Williams (1987), Theorem IV.45.9 on page 105] applied to $f(x) = x^{p/q}$, for $t \leq \tau$, we have

$$dY^p = d((Y^q)^{p/q}) = \frac{p}{q}Y^{p-q}\,dY^q + \frac{1}{2}\frac{p}{q}\left(\frac{p}{q}-1\right)Y^{p-2q}I_{Y>0}\,d[Y^q]$$
$$= \frac{p}{q}Y^{p-q}\left[-2qY^{q-1/2}\gamma\sum_j X_j\,dB_j\right.$$
$$\left. + 2q|X|^2Y^{q-1}I_{Y>0}[g + (q-1)\gamma^2]\,dt - qn\gamma^2 Y^q\,dt + d\varphi_t^{(q)}\right]$$
$$+ \frac{1}{2}\frac{p}{q}\left(\frac{p}{q}-1\right)Y^{p-2q}I_{Y>0}[4q^2 Y^{2q-1}\gamma^2|X|^2]\,dt$$

$$= -2pY^{p-1/2}\gamma \sum_j X_j \, dB_j$$
$$+ 2p|X|^2 Y^{p-1} I_{Y>0}[g + (q-1)\gamma^2] \, dt$$
$$- pn\gamma^2 Y^p \, dt + 0$$
$$+ 2p(p-q)Y^{p-1} I_{Y>0} \gamma^2 |X|^2 \, dt$$
$$= dY^p - d\varphi^{(p)} \qquad [\text{by } (2.8)].$$

Thus, $\varphi_t^{(p)} \equiv 0$, as claimed. $\square$

The last result of this section is needed in the proof of Theorem 1.1.

LEMMA 2.2. *Any solution $X$ of (1.2) spends zero Lebesgue time on the boundary:*

$$\int_0^t I_{|X_s|=1} \, ds = 0 \qquad a.s.$$

PROOF. With $\tau$ from Lemma 2.1, it suffices to show

$$\int_0^\tau I_{|X_s|=1} \, ds = 0.$$

Applying the extension of Itô's formula to $C^1$ functions (cited above) to $f(x) = x^{1/p}$, for $t \leq \tau$, using Lemma 2.1,

(2.10)
$$dY = d((Y^p)^{1/p}) = \frac{1}{p} Y^{1-p} \, dY^p + \frac{1}{2}\frac{1}{p}\left(\frac{1}{p} - 1\right) Y^{1-2p} I_{Y>0} \, d[Y^p]$$
$$= -2\gamma Y^{1/2} \sum_j X_j \, dB_j$$
$$+ 2|X|^2 I_{Y>0}[g + (p-1)\gamma^2] \, dt$$
$$- n\gamma^2 Y \, dt$$
$$+ 2(1-p)\gamma^2 I_{Y>0} |X|^2 \, dt$$
$$= -2\gamma Y^{1/2} \sum_j X_j \, dB_j + 2|X|^2 I_{Y>0} g \, dt - n\gamma^2 Y \, dt.$$

On the other hand,

$$dY = d[1 - |X|^2]$$
$$= -2Y^{1/2}\gamma \sum_j X_j \, dB_j + 2g|X|^2 \, dt - nY\gamma^2 \, dt.$$

Upon comparison with (2.10), we must have
$$2g|X|^2 I_{Y=0}\, dt = 0,$$
which is equivalent to
$$2g(1) I_{|X|=1}\, dt = 0,$$
since $Y = 0 \Leftrightarrow |X| = 1$. The desired conclusion follows because $g(1) > 0$. □

**3. Proof of Theorem 1.1.** Let $X$, $\widetilde{X}$ be solutions of (1.2) with the same underlying Brownian motion. Suppose $p \in (\frac{1}{2}, 1)$ satisfies $p > 1 - \frac{g(1)}{\gamma^2(1)}$ and choose $\varepsilon = \varepsilon(p)$ as in (2.1). It is no loss to assume $\varepsilon < \frac{1}{2}$. Let $Y = 1 - |X|^2$, $\widetilde{Y} = 1 - |\widetilde{X}|^2$ and

$$W = (Y^p - \widetilde{Y}^p)^2 + |X - \widetilde{X}|^2, \tag{3.1}$$

as in the Introduction.

Since the coefficients of (1.2) are locally Lipschitz on the interior of the unit ball $E$, pathwise uniqueness holds up to the first hit of the boundary. Hence, it suffices by a restart argument to consider starting points on the boundary. By making $\varepsilon$ smaller if necessary, it suffices to show for $\tau = \min(\tau_\varepsilon(Y), \tau_\varepsilon(\widetilde{Y}))$,
$$X_t = \widetilde{X}_t, \qquad t \leq \tau,\ \text{a.s.}$$

Below we will use the fact
$$Y_t \vee \widetilde{Y}_t \leq \varepsilon, \qquad t \leq \tau. \tag{3.2}$$

Writing
$$G(u) = g(u) + (p-1)\gamma^2(u),$$
by Lemma 2.1,

$$
\begin{aligned}
d(Y^p - \widetilde{Y}^p) &= -2p \sum_j [\gamma(|X|) Y^{p-1/2} X_j - \gamma(|\widetilde{X}|) \widetilde{Y}^{p-1/2} \widetilde{X}_j]\, dB_j \\
&\quad + 2p[|X|^2 Y^{p-1} I_{Y>0} G(|X|) - |\widetilde{X}|^2 \widetilde{Y}^{p-1} I_{\widetilde{Y}>0} G(|\widetilde{X}|)]\, dt \\
&\quad - np[\gamma^2(|X|) Y^p - \gamma^2(|\widetilde{X}|) \widetilde{Y}^p]\, dt \\
&= dM + I_1\, dt + I_2\, dt \qquad \text{say.}
\end{aligned}
\tag{3.3}
$$

Then
$$
\begin{aligned}
d(Y^p - \widetilde{Y}^p)^2 &= 2(Y^p - \widetilde{Y}^p)[dM + I_1\, dt + I_2\, dt] \\
&\quad + 4p^2 \sum_j [\gamma(|X|) Y^{p-1/2} X_j - \gamma(|\widetilde{X}|) \widetilde{Y}^{p-1/2} \widetilde{X}_j]^2\, dt \\
&= 2(Y^p - \widetilde{Y}^p)[dM + I_1\, dt + I_2\, dt] + I_3\, dt \qquad \text{say.}
\end{aligned}
\tag{3.4}
$$

The term $(Y^p - \widetilde{Y}^p)I_1 \, dt$ is the "good" negative singular term that will compensate for the "bad" positive singular term $I_3 \, dt$. Note $I_3$ is "singular" in the sense that $I_3/W$ is unbounded. It turns out the term involving $I_2$ is not singular in this sense. In Swart's article (i.e., $p = \frac{1}{2}$) the term involving $I_3$ is not singular in that $I_3/W$ is bounded in this case.

We also have

$$
\begin{aligned}
d|X - \widetilde{X}|^2 &= 2[Y^{1/2}\gamma(|X|) - \widetilde{Y}^{1/2}\gamma(|\widetilde{X}|)] \sum_j (X_j - \widetilde{X}_j) \, dB_j \\
&\quad - 2\sum_j (X_j - \widetilde{X}_j)[g(|X|)X_j - g(|\widetilde{X}|)\widetilde{X}_j] \, dt \\
&\quad + n[Y^{1/2}\gamma(|X|) - \widetilde{Y}^{1/2}\gamma(|\widetilde{X}|)]^2 \, dt \\
&= dR + I_4 \, dt + I_5 \, dt \quad \text{say.}
\end{aligned}
\tag{3.5}
$$

Exactly as in Swart, the term $I_4 \, dt$ is nonsingular and the term $I_5 \, dt$ is positive and singular.

Now we estimate the individual terms. For notational convenience write

$$
Z = (Y^p - \widetilde{Y}^p)(\widetilde{Y}^{p-1} - Y^{p-1}).
\tag{3.6}
$$

Note $Z \geq 0$ and it will turn out all the singular terms involve $Z$.

LEMMA 3.1. *For $Y$ and $\widetilde{Y}$ positive and $t \leq \tau$,*

$$(Y^p - \widetilde{Y}^p)I_1 \leq -2pZ|X|^2 G(|X|) + C\varepsilon Z,$$

*where $C$ is independent of $\varepsilon$.*

PROOF. By (3.3),

$$
\begin{aligned}
(Y^p - \widetilde{Y}^p)I_1 &= (Y^p - \widetilde{Y}^p)2p[|X|^2 Y^{p-1}G(|X|) - |\widetilde{X}|^2 \widetilde{Y}^{p-1}G(|\widetilde{X}|)] \\
&= 2p(Y^p - \widetilde{Y}^p)\{(Y^{p-1} - \widetilde{Y}^{p-1})|X|^2 G(|X|) \\
&\quad + \widetilde{Y}^{p-1}[|X|^2 G(|X|) - |\widetilde{X}|^2 G(|\widetilde{X}|)]\} \\
&= -2pZ|X|^2 G(|X|) \\
&\quad + 2p(Y^p - \widetilde{Y}^p)\widetilde{Y}^{p-1}[|X|^2 G(|X|) - |\widetilde{X}|^2 G(|\widetilde{X}|)].
\end{aligned}
\tag{3.7}
$$

Thus, we need only estimate the last term.

Recall we are assuming $\varepsilon < \frac{1}{2}$. Then, for $t \leq \tau$, by (3.2),

$$|X_t| = \sqrt{1 - Y_t} \geq \sqrt{1 - \varepsilon} \geq \sqrt{\tfrac{1}{2}}$$

and this holds for $|\widetilde{X}_t|$ too. Thus,

$$
\begin{aligned}
||X| - |\widetilde{X}|| &= \frac{||X|^2 - |\widetilde{X}|^2|}{|X| + |\widetilde{X}|} \\
&\leq \tfrac{\sqrt{2}}{2}||X|^2 - |\widetilde{X}|^2| \\
&= \tfrac{\sqrt{2}}{2}|Y - \widetilde{Y}|.
\end{aligned}
\tag{3.8}
$$

Now $u^2 G(u)$ is Lipschitz, hence, for some constant $C$ independent of $\varepsilon$, the last term in (3.7) is bounded by

$$C|Y^p - \widetilde{Y}^p|\widetilde{Y}^{p-1}||X| - |\widetilde{X}|| \leq C|Y^p - \widetilde{Y}^p|\widetilde{Y}^{p-1}\tfrac{\sqrt{2}}{2}|Y - \widetilde{Y}|.$$

Since $Y \vee \widetilde{Y} \leq \varepsilon$, we just need to show for some $C > 0$ independent of $\varepsilon$,

$$\widetilde{Y}^{p-1}|Y^p - \widetilde{Y}^p||Y - \widetilde{Y}| \leq C\max(Y,\widetilde{Y})|Y^p - \widetilde{Y}^p||Y^{p-1} - \widetilde{Y}^{p-1}|$$
$$[= C\max(Y,\widetilde{Y})Z].$$

To this end, since $p - 1 < 0$, the worst case occurs if $\widetilde{Y} \leq Y$. Thus, it is enough to show

$$\widetilde{Y}^{p-1}(Y - \widetilde{Y}) \leq CY(\widetilde{Y}^{p-1} - Y^{p-1}),$$

which after division by $\widetilde{Y}^{p-1}Y$ is equivalent to

$$1 - \frac{\widetilde{Y}}{Y} \leq C\left[1 - \left(\frac{\widetilde{Y}}{Y}\right)^{1-p}\right].$$

It is easy to see $\sup_{0 < u < 1} \frac{1-u}{1-u^{1-p}} < \infty$, and so taking $C$ to be the supremum does the trick. □

LEMMA 3.2. *For $Y$ and $\widetilde{Y}$ positive and $t \leq \tau$,*

$$|I_2| \leq C[|Y^p - \widetilde{Y}^p| + |X - \widetilde{X}|],$$

*where $C$ is independent of $\varepsilon$.*

PROOF. From (3.3),

$$|I_2| \leq np[|Y^p - \widetilde{Y}^p|\gamma^2(|X|) + \widetilde{Y}^p|\gamma^2(|X|) - \gamma^2(|\widetilde{X}|)|].$$

Since $\widetilde{Y} < \varepsilon < \tfrac{1}{2}$ and $\gamma$ is bounded and Lipschitz, the latter is bounded by

$$np|Y^p - \widetilde{Y}^p|\sup\gamma^2 + C||X| - |\widetilde{X}||,$$

where $C$ is independent of $\varepsilon$. □

In order to bound the singular term $I_3$, we need the following lemma.

LEMMA 3.3. *Let $x, y > 0$, $p \in (\frac{1}{2}, 1)$ and set $z = (x^p - y^p)(y^{p-1} - x^{p-1})$. Then, for some constant $C$ depending only on $p$,*

$$|x^{p-1/2} - y^{p-1/2}||x - y| \leq C \max(x^{1/2} y^{1-p}, y^{1/2} x^{1-p}) z.$$

PROOF. It is no loss to assume $y < x$. Then

$$\frac{(x^{p-1/2} - y^{p-1/2})(x-y)}{z} = \frac{(x^{p-1/2} - y^{p-1/2})(x-y) x^{1-p} y^{1-p}}{(x^p - y^p)(x^{1-p} - y^{1-p})}$$

$$= x^{1/2} y^{1-p} \frac{(1 - (y/x)^{p-1/2})(1 - y/x)}{(1 - (y/x)^p)(1 - (y/x)^{1-p})}$$

$$\leq x^{1/2} y^{1-p} \sup_{0 < u < 1} \frac{(1 - w^{p-1/2})(1 - w)}{(1 - w^p)(1 - w^{1-p})}.$$

The supremum is finite because $p \in (\frac{1}{2}, 1)$. Taking $C$ as this value yields the desired conclusion. □

LEMMA 3.4. *For $Y$ and $\widetilde{Y}$ positive and $t \leq \tau$,*

$$I_3 \leq \frac{p(2p-1)^2}{(1-p)} \gamma^2(|X|) |X|^2 Z + C|X - \widetilde{X}|^2 + C\varepsilon Z,$$

*where $C$ is independent of $\varepsilon$.*

PROOF. To ease eye strain, write

$$a = Y^{p-1/2},$$
$$b = \widetilde{Y}^{p-1/2},$$
$$U = \gamma(|X|) X,$$
$$V = \gamma(|\widetilde{X}|) \widetilde{X}.$$

Then from (3.4),

$$I_3 = 4p^2 |aU - bV|^2$$

and $a, b \in (0, 1)$. Since $ab < 1$,

$$I_3 - 4p^2 |U - V|^2$$

$$= 4p^2 \left[ (a^2 - 1)|U|^2 - 2(ab - 1) \sum_i U_i V_i + (b^2 - 1)|V|^2 \right]$$

$$\leq 4p^2 [(a^2 - 1)|U|^2 - 2(ab - 1)|U||V| + (b^2 - 1)|V|^2]$$

(3.9) $\qquad = 4p^2 [(a|U| - b|V|)^2 - (|U| - |V|)^2]$

$$\leq 4p^2[a|U| - b|V|]^2$$
$$= 4p^2[(a-b)|U| + b(|U| - |V|)]^2$$
$$= 4p^2[(a-b)^2|U|^2 + 2b(a-b)|U|(|U| - |V|) + b^2(|U| - |V|)^2].$$

Since $r\gamma(r)$ is Lipschitz in $r$, for some $C > 0$ independent of $\varepsilon$,

$$||U| - |V|| \leq C||X| - |\widetilde{X}||.$$

Hence, the last term on left-hand side of (3.9) is bounded by $4p^2 C|X - \widetilde{X}|^2$ and the middle term is bounded by

$$8p^2 Cb|a - b|||X| - |\widetilde{X}||$$
$$= 8p^2 C\widetilde{Y}^{p-1/2}|Y^{p-1/2} - \widetilde{Y}^{p-1/2}|||X| - |\widetilde{X}||$$
$$\leq C\widetilde{Y}^{p-1/2}|Y^{p-1/2} - \widetilde{Y}^{p-1/2}||Y - \widetilde{Y}| \quad \text{[by (3.8)]}$$
$$\leq C\varepsilon^{p-1/2}\varepsilon^{1/2}\varepsilon^{1-p}Z \quad \text{[by Lemma 3.3 and (3.6)]}$$
$$= C\varepsilon Z.$$

Also, since $\gamma$ is bounded and Lipschitz,

$$|U - V| \leq |\gamma(|X|)(X - \widetilde{X})| + |\widetilde{X}(\gamma(|X|) - \tilde{\gamma}(|X|))|$$
$$\leq C[|X - \widetilde{X}| + ||X| - |\widetilde{X}||]$$
$$\leq C|X - \widetilde{X}|.$$

Thus, to finish the proof we just need to show the first term on right-hand side of (3.9) satisfies

$$4p^2(a-b)^2|U|^2 \leq \frac{p(2p-1)^2}{1-p}\gamma^2(|X|)|X|^2 Z.$$

But, by Lemma A.1,

$$(a-b)^2 = (Y^{p-1/2} - \widetilde{Y}^{p-1/2})^2$$
$$\leq \frac{(2p-1)^2}{4p(1-p)}Z$$

and since $|U|^2 = \gamma^2(|X|)|X|^2$, we get the desired bound. $\square$

LEMMA 3.5. *For $Y$ and $\widetilde{Y}$ positive and $t \leq \tau$,*

$$|I_4| \leq C|X - \widetilde{X}|^2,$$

*where $C$ is independent of $\varepsilon$.*

PROOF. Since $g$ is Lipschitz, the bound is clear from the definition (3.5) of $I_4$. $\square$

Finally, we bound the last singular term $I_5$.

LEMMA 3.6. *For positive $Y$ and $\widetilde{Y}$ and $t \leq \tau$,*
$$I_5 \leq C[\varepsilon^{2-2p} Z + |X - \widetilde{X}|^2],$$
*where $C$ is independent of $\varepsilon$.*

PROOF. By (3.5),
$$I_5 = n[(Y^{1/2} - \widetilde{Y}^{1/2})\gamma(|X|) + \widetilde{Y}^{1/2}(\gamma(|X|) - \gamma(|\widetilde{X}|))]^2$$
$$\leq 2n(\sup \gamma^2)(Y^{1/2} - \widetilde{Y}^{1/2})^2 + 2n\widetilde{Y}(\gamma(|X|) - \gamma(|\widetilde{X}|))^2.$$

Since $\gamma$ is Lipschitz and $\widetilde{Y} \leq 1$, for some $C > 0$ independent of $\varepsilon$, the second term is bounded by
$$C|X - \widetilde{X}|^2.$$

Hence, we need only show
$$(Y^{1/2} - \widetilde{Y}^{1/2})^2 \leq C \varepsilon^{2-2p}(Y^p - \widetilde{Y}^p)(\widetilde{Y}^{p-1} - Y^{p-1})$$

(since the latter is $C\varepsilon^{2-2p} Z$) where $C$ is independent of $\varepsilon$. To this end, write $x = Y^{1/2}$ and $y = \widetilde{Y}^{1/2}$ and without loss of generality, assume $y < x$. Then since $x, y < \sqrt{\varepsilon}$ [by (3.2)], it is enough to show
$$(x - y)^2 \leq C x^{2-2p} y^{2-2p}(x^{2p} - y^{2p})(y^{2p-2} - x^{2p-2})$$

or, equivalently,
$$(x - y)^2 \leq C(x^{2p} - y^{2p})(x^{2-2p} - y^{2-2p}).$$

By dividing by $x^2$, this is equivalent to
$$\left(1 - \frac{y}{x}\right)^2 \leq C\left(1 - \left(\frac{y}{x}\right)^{2p}\right)\left(1 - \left(\frac{y}{x}\right)^{2-2p}\right).$$

Since $p < 1$, it is easy to see
$$C := \sup_{0 < w < 1} \frac{(1-w)^2}{(1 - w^{2p})(1 - w^{2-2p})} < \infty$$

does the trick. $\square$

Now we can show how the hypothesis $\frac{g(1)}{\gamma^2(1)} > \sqrt{2} - 1$ implies the negative singular term $(Y^p - \widetilde{Y}^p) I_1$ compensates for the positive singular terms $I_3$ and $I_5$.

LEMMA 3.7. *Let $Y$ and $\widetilde{Y}$ be positive and suppose $\frac{g(1)}{\gamma^2(1)} > \sqrt{2} - 1$. Then, for $p = 1 - \frac{\sqrt{2}}{4}$, by making $\varepsilon$ smaller if necessary, for $t \leq \tau$, we have*

$$2(Y^p - \widetilde{Y}^p)I_1 + I_3 + I_5 \leq C|X - \widetilde{X}|^2,$$

*where $C > 0$ is independent of $\varepsilon$.*

PROOF. First note the value $p = 1 - \frac{\sqrt{2}}{4} \in (\frac{1}{2}, 1)$ minimizes the function

$$-(p-1) + \frac{(2p-1)^2}{4(1-p)}$$

on $(\frac{1}{2}, 1)$ and the minimum value is $\sqrt{2} - 1$. As shown below, this is what leads to the hypothesis $\frac{g(1)}{\gamma^2(1)} > \sqrt{2} - 1$. Note too that this choice of $p$ satisfies $p > 1 - \frac{g(1)}{\gamma(1)^2}$ (so Lemma 2.1 is applicable).

By Lemmas 3.1, 3.4 and 3.6,

$$\begin{aligned}
2(Y^p - \widetilde{Y}^p)I_1 + I_3 + I_5 \\
\leq -4pZ|X|^2 G(|X|) + C\varepsilon Z \\
+ \frac{p(2p-1)^2}{1-p}\gamma^2(|X|)|X|^2 Z + C|X - \widetilde{X}|^2 + C\varepsilon Z
\end{aligned}$$

(3.10)

$$\begin{aligned}
+ C[\varepsilon^{2-2p}Z + |X - \widetilde{X}|^2] \\
= 4pZ|X|^2\left[-G(|X|) + \frac{(2p-1)^2}{4(1-p)}\gamma^2(|X|)\right] \\
+ C[\varepsilon + \varepsilon^{2-p}]Z + C|X - \widetilde{X}|^2.
\end{aligned}$$

To finish, we show the coefficient of $Z$ is negative for $\varepsilon$ sufficiently small. Recall $G(u) = g(u) + (p-1)\gamma^2(u)$, hence,

$$\begin{aligned}
-G(u) + \frac{(2p-1)^2}{4(1-p)}\gamma^2(u) &= \gamma^2(u)\left[-\frac{g(u)}{\gamma^2(u)} - (p-1) + \frac{(2p-1)^2}{4(1-p)}\right] \\
&= \gamma^2(u)\left[-\frac{g(u)}{\gamma^2(u)} + \sqrt{2} - 1\right] \qquad \text{by choice of } p.
\end{aligned}$$

This is where we see the reason for choosing $p$ to minimize $-(p-1) + \frac{(2p-1)^2}{4(1-p)}$. Then the coefficient $K$ of $Z$ in (3.10) is

$$K = 4p|X|^2\gamma^2(|X|)\left[-\frac{g(|X|)}{\gamma^2(|X|)} + \sqrt{2} - 1\right] + C(\varepsilon + \varepsilon^{2-p}),$$

which by (3.2) (and that $|X| = \sqrt{1-Y}$) is bounded by

$$4p|X|^2\gamma^2(|X|)\left[-\inf_{u^2 \geq 1-\varepsilon} \frac{g(u)}{\gamma^2(u)} + \sqrt{2} - 1\right] + C(\varepsilon + \varepsilon^{2-p}).$$

By our hypotheses $\frac{g(1)}{\gamma^2(1)} > \sqrt{2} - 1$ and continuity of $g$ and $\gamma$, and by making $\varepsilon$ smaller if necessary, the quantity in square brackets is negative, hence,

$$K \leq 4p\ (1-\varepsilon)[\inf \gamma^2]\bigg[-\inf_{u^2 \geq 1-\varepsilon} \frac{g(u)}{\gamma^2(u)} + \sqrt{2} - 1\bigg] + C(\varepsilon + \varepsilon^{2-p}).$$

Then for the same reason, making $\varepsilon$ smaller if necessary, $K < 0$. $\square$

Now we prove Theorem 1.1. With $W$ from (3.1) and $p$ from Lemma 3.7, for $Y$ and $\widetilde{Y}$ positive, by (3.4) and (3.5),

$$dW = 2(Y^p - \widetilde{Y}^p)[dM + I_1\,dt + I_2\,dt] + I_3\,dt$$
$$+ dR + I_4\,dt + I_5\,dt.$$

Hence, upon integrating and using Lemma 2.2, we can apply Lemmas 3.2, 3.5 and 3.7 to get

$$E[W_{t\wedge\tau}] \leq CE\bigg[\int_0^{t\wedge\tau} \{|X_s - \widetilde{X}_s|^2 + (Y_s^p - \widetilde{Y}_s^p)^2 + |Y_s^p - \widetilde{Y}_s^p||X_s - \widetilde{X}_s|\}\,ds\bigg]$$
$$\leq 2CE\bigg[\int_0^{t\wedge\tau} W_s\,ds\bigg].$$

Then Gronwall's inequality yields $W_{t\wedge\tau} = 0$ a.s., which forces $X_{t\wedge\tau} = \widetilde{X}_{t\wedge\tau}$ a.s., as desired.

REMARK 3.1. We have written the proof of Theorem 1.1 in such a way that it is easy to change to the situation of Theorem 1.3. One uses $Y = h(X)$ in place of $1 - |X|^2$ in the argument and computes the differential of powers of $h(X)$ to obtain the analog of Lemma 2.1. Again, our proofs in Section 2 are given with this in mind. The notation is more complex, but the basic ideas are the same.

REMARK 3.2. If one replaces the power $\frac{1}{2}$ in (1.2) by $r \in (\frac{1}{2}, 1)$, then, as pointed out in the Introduction, our method can be used to show pathwise uniqueness holds for starting points on the boundary (recall the boundary is unattainable for all other starting points and pathwise uniqueness follows from standard results). In this case there will be no restriction such as $g(1)/\gamma^2(1) > \sqrt{2} - 1$. This is due to the unattainable nature of the boundary previously mentioned. The analog of Lemma 2.1 exemplifies this: one can take any $p \in (1-r, 1)$. Consequently, the proofs of the analogs of Lemmas 3.4 and 3.7 are simpler. The remaining details are similar to those furnished above.

REMARK 3.3. There is the question of *existence* of a solution to the equation (1.4) considered in Theorem 1.3. Basically, it is necessary to verify the positive maximum principle holds for the operator

$$\frac{1}{2}\sum_{i,j} h^{1/2} a_{ij} \frac{\partial^2}{\partial x_i \, \partial x_j} + \sum_i b_i \frac{\partial}{\partial x_i}.$$

See Theorem 4.5.4 (page 199) in Ethier and Kurtz (1986) where this is done for the martingale problem. By problem 19 (page 265) from the same chapter and Theorem 5.3.3 (page 293), this can be translated into existence of weak solutions to SDEs. This argument is due to the referee.

**4. Proof of Theorem 1.2.** We prove the equivalent statement that weak uniqueness holds for the SDE

$$dX_t = [1 - |X_t|^2]^r \gamma(1 - |X_t|^2) \, dB_t - g(1 - |X_t|^2) X_t \, dt,$$
$$X_0 = x$$

in $E$, where $\gamma$ and $g$ are continuous and strictly positive on $[0,1]$. Since the diffusion matrix is uniformly positive definite on compact subsets contained in the interior of $E$, it is well known that weak uniqueness holds up to the first hitting time of the boundary. Thus, it suffices to show weak uniqueness holds for starting points on the boundary, and by rotational invariance, it is no loss to take

$$X_0 = (0, \ldots, 0, 1).$$

First, transform the state space. The mapping

(4.1)
$$v = v(x) = 1 - |x|^2,$$
$$y = y(x) = \left(\frac{x_1}{|x|}, \ldots, \frac{x_{n-1}}{|x|}\right), \qquad x \in B_1(1) \cap E,$$

is one-to-one. We want to compute the differentials of

(4.2)
$$V_t = v(X_t),$$
$$Y_t = y(X_t).$$

To this end, for $y \in \mathbb{R}^{n-1}$ and $i,j \in \{1, \ldots, n-1\}$, define

(4.3)
$$A_{ij}(y) = \begin{cases} 1 - y_i^2, & i = j, \\ -y_i y_j, & i \neq j. \end{cases}$$

LEMMA 4.1. *For $|y| \leq \frac{1}{2}$, $A(y)$ is uniformly positive definite: for some constant $\lambda > 0$,*

$$\langle A(y)\xi, \xi \rangle \geq \lambda |\xi|^2, \qquad \xi \in \mathbb{R}^{n-1}.$$

PROOF. For $y, \xi \in \mathbb{R}^{n-1}$ with $|y| \leq \frac{1}{2}$,

$$\langle A(y)\xi, \xi \rangle = \sum_{i,j=1}^{n-1} A_{ij}(y)\xi_i\xi_j$$

$$= \sum_i (1 - y_i^2)\xi_i^2 - 2\sum_{i<j} y_i y_j \xi_i \xi_j$$

$$= |\xi|^2 - \langle \xi, y \rangle^2$$

$$\geq |\xi|^2 - |\xi|^2 |y|^2$$

$$\geq \tfrac{3}{4}|\xi|^2. \qquad \square$$

Using the formulas

$$\frac{\partial y_i}{\partial x_j} = \begin{cases} \dfrac{|x|^2 - x_i^2}{|x|^3}, & i = j, \\ -\dfrac{x_i x_j}{|x|^3}, & i \neq j, \end{cases}$$

$$\frac{\partial^2 y_i}{\partial x_j^2} = \begin{cases} -\dfrac{3x_i(|x|^2 - x_i^2)}{|x|^5}, & i = j, \\ -\dfrac{x_i(|x|^2 - 3x_j^2)}{|x|^5}, & i \neq j, \end{cases}$$

Itô's formula can be used to show that

$$dV = -2V^r \gamma(V)\sqrt{1-V}\,d\beta + [2g(V)(1-V) - nV^{2r}\gamma^2(V)]\,dt$$

(4.4) $\quad dY = V^r \gamma(V)(1-V)^{-1/2} A^{1/2}(Y)\,dM$

$$- \frac{n-1}{2}(1-V)^{-1}V^{2r}\gamma^2(V)Y\,dt,$$

where $(\beta, M) \in \mathbb{R} \times \mathbb{R}^{n-1}$ is an $n$-dimensional Brownian motion. Hence, it suffices to show weak uniqueness holds for (4.4), with $(V_0, Y_0) = (0, 0)$.

In what follows, we use the last chapter in Breiman (1968) as our basic reference for one-dimensional diffusions. Notice $V$ is an autonomous one-dimensional diffusion process with state space $[0, 1]$. Since the process $X$ never hits 0 when $X_0 \neq 0$ (since $n \geq 2$), the state space will actually be $[0, 1)$ since we are taking $V_0 = 0$. Transform the state space using the scale function $s(v)$ given by $s(0) = 0$ and

$$s'(v) = \exp\left(-\int_0^v \frac{2g(w)(1-w) - nw^{2r}\gamma^2(w)}{2w^{2r}\gamma^2(w)(1-w)}\,dw\right), \qquad v \in [0, 1).$$

Since $r < \frac{1}{2}$, the integral $\int_0^v$ is finite (again, using $n \geq 2$) and, therefore, it turns out $s(V)$ is a diffusion in the natural scale with state space $[0, \infty)$, and 0 is a slowly reflecting boundary point (i.e., the process spends positive Lebesgue time there). By a time change we can convert $s(V)$ into a one-dimensional Brownian motion in $[0, \infty)$ with slow reflection at 0.

Using this time change on $(s(V), Y)$, weak uniqueness for (4.4), with $(V_0, Y_0) = (0, 0)$, comes down to proving weak uniqueness for

$$
\begin{aligned}
dU &= d\beta + \tfrac{1}{2}\ell_t^0(U), \\
dN &= \frac{1}{2}A^{1/2}(N)H(U)I_{U>0}\,dM - \frac{n-1}{8}H^2(U)NI_{U>0}\,dt, \\
I_{U=0}\,dt &= c\,d\ell_t^0(U), \\
(U_0, N_0) &= (0, 0)
\end{aligned}
\tag{4.5}
$$

in $\mathbb{R}_+ \times \mathbb{R}^{n-1}$, where $s^{-1}$ is the inverse of $s$,

$$H(U) = [s' \circ s^{-1}(U)]^{-1}(1 - s^{-1}(U))^{-1},$$

$c > 0$ and $\ell_t^0(U)$ is the local time of $U$ at 0.

To prove that weak uniqueness holds for (4.5) we introduce a certain stopped submartingale problem. Let $\Omega = C([0, \infty), \mathbb{R}^n)$ be the space of continuous paths in $\mathbb{R}^n$ and equip it with the cylindrical Borel $\sigma$-algebra. Denote by $Z_t(\omega)$ the coordinate process $Z_t(\omega) = \omega(t)$, $\omega \in \Omega$, and let $\mathcal{F}_t = \sigma(Z_s : s \leq t)$, $\mathcal{F} = \sigma(Z_s : s \geq 0)$. For $(u, y) \in \mathbb{R}^n$ with $y = (y_1, \ldots, y_{n-1})$, set

$$L = \frac{1}{2}\left[\frac{\partial^2}{\partial u^2} + \frac{1}{4}\sum_{i,j=1}^{n-1} A_{ij}(y)H^2(U)\frac{\partial^2}{\partial y_i\,\partial y_j}\right] - \frac{n-1}{8}H^2(U)\sum_{i=1}^{n-1} y_i\frac{\partial}{\partial y_i},$$

where $A$ is from (4.3). A probability measure $P$ on $(\Omega, \mathcal{F})$ solves the *stopped submartingale problem* if for the first exit time $\tau$ of $Z$ from a small neighborhood of $(0, 0)$ in $\mathbb{R}_+ \times \mathbb{R}^{n-1}$, we have

$$P(Z_0 = 0) = 1,$$
$$P(Z_{t \wedge \tau} \in \mathbb{R}_+ \times \mathbb{R}^{n-1}) = 1$$

and for all $f \in C^{1,2}([0, \infty) \times \mathbb{R}^n)$ satisfying

$$c\frac{\partial f}{\partial t} + \frac{1}{2}\frac{\partial f}{\partial u} \geq 0 \qquad \text{on } [0, \infty) \times \{(u, y) : u = 0, y \in \mathbb{R}^{n-1}\}, \tag{4.6}$$

we have

$$f(t \wedge \tau, Z_{t \wedge \tau}) - \int_0^{t \wedge \tau} I_{U_s > 0}\left[\frac{\partial f}{\partial s} + Lf\right](s, Z_s)\,ds$$

is a $P$-submartingale.

The matrix of coefficients of second-order terms in $L$ is bounded, continuous and uniformly elliptic in a neighborhood of $(0, 0)$ in $\mathbb{R}_+ \times \mathbb{R}^{n-1}$ (using

Lemma 4.1), and the first-order term coefficients are continuous there. The boundary operator in (4.6) has the form

$$\rho \frac{\partial}{\partial t} + \gamma \cdot \nabla,$$

where $\rho$ is continuous and positive and $\gamma$ is Lipschitz continuous. Moreover, if $n$ is the unit inward normal to $\partial(\mathbb{R}_+ \times \mathbb{R}^{n-1})$, then $\gamma \cdot n > 0$. Thus, by Theorem 5.8 on page 196 of Stroock and Varadhan (1971), uniqueness holds for the stopped submartingale problem.

A routine use of Itô's formula shows any solution to (4.5) yields a solution of the stopped submartingale problem. Hence, weak uniqueness holds for (4.5), as desired.

## 5. Open questions.

1. Our method slightly improved Swart's condition $c \geq 1$ for (1.1), but is still remains to resolve the case $0 < c \leq 2(\sqrt{2} - 1)$ for $n \geq 2$. This seems quite difficult.
2. Is the assumption in Theorem 1.3 that $g|\nabla h|$ and $\langle a \nabla h, \nabla h \rangle$ are Lipschitz functions of $h$ really needed? Our proof breaks down without it.
3. The question of pathwise uniqueness for

$$dX_t = [1 - |X_t|^2]^r \gamma(|X_t|) \, dB_t - g(|X_t|) X_t \, dt,$$

with $r < \frac{1}{2}$ and Lipschitz $\gamma, g : [0,1] \to (0,1)$, remains open.
4. With reference to the equation studied by Athreya, Barlow, Bass and Perkins (2002) and Bass and Perkins (2002) mentioned in the Introduction, decide whether or not pathwise uniqueness holds for

$$dX_t^j = \sum_k (2X_t^i)^{1/2} \sigma_{ik}(X_t) \, dB_t^k + b^i(X_t) \, dt, \qquad i = 1, \ldots, n,$$

in the positive orthant.

## APPENDIX: A TECHNICAL RESULT

LEMMA A.1. *For positive $x$ and $y$ and $p \in (\frac{1}{2}, 1)$,*

$$[x^{p-1/2} - y^{p-1/2}]^2 \leq \frac{(2p-1)^2}{4p(1-p)} (x^p - y^p)(y^{p-1} - x^{p-1}).$$

PROOF. It is no loss to assume $y < x$. Then

$$\frac{(x^{p-1/2} - y^{p-1/2})^2}{(x^p - y^p)(y^{p-1} - x^{p-1})} = \frac{(y/x)^{p-1}(1 - (y/x)^{p-1/2})^2}{(1 - (y/x)^p)(1 - (y/x)^{1-p})}.$$

Thus, we must show

$$\sup_{0<w<1} \frac{w^{1-p}(1-w^{p-1/2})^2}{(1-w^p)(1-w^{1-p})} = \frac{(2p-1)^2}{4p(1-p)}$$

and replacing $w$ by $z^2$ and then $2p$ by $q$, this comes down to showing

(A.1) $$\sup_{0<z<1} \frac{z^{2-q}(1-z^{q-1})^2}{(1-z^q)(1-z^{2-q})} = \frac{(q-1)^2}{q(2-q)}, \qquad 1 < q < 2.$$

If we show

$$f(z) = \frac{(1-z^q)(1-z^{2-q})}{z^{2-q}(1-z^{q-1})^2}, \qquad 0 < z < 1,$$

is decreasing, then since

$$\lim_{z \to 1^-} f(z) = \frac{q(2-q)}{(q-1)^2},$$

(A.1) will hold. The proof is elementary, but a bit involved. Define

$$f_1(z) = (q+1)z - q - z^{1-q}, \qquad z > 0.$$

Since $f_1'(z) = q + 1 + (q-1)z^{-q} > 0$, we have

(A.2) $$f_1(z) \leq f_1(1) = 0, \qquad z \in (0,1].$$

If

$$f_2(z) = -(q+1)(2-q)z^q - q(q-1)z^{q-1} + q(3-q)z - (q-1)(2-q),$$
$$z \geq 0,$$

then

$$f_2'(z) = qz^{q-2}[-(q+1)(2-q)z - (q-1)^2 + (3-q)z^{2-q}], \qquad z > 0.$$

If there is $z_0 \in (0,1)$ such that $f_2'(z_0) = 0$, then

$$-(q+1)(2-q)z_0 - (q-1)^2 + (3-q)z_0^{2-q} = 0$$

and this yields

$$(3-q)z_0^{2-q} = (q+1)(2-q)z_0 + (q-1)^2.$$

Hence,

$$f_2(z_0) = -(q+1)(2-q)z_0^q - q(q-1)z_0^{q-1} + q(3-q)z_0 - (q-1)(2-q)$$
$$= z_0^{q-1}[-(q+1)(2-q)z_0 - q(q-1)$$
$$\qquad + q(3-q)z_0^{2-q} - (q-1)(2-q)z_0^{1-q}]$$

$$
\begin{aligned}
&= z_0^{q-1}[-(q+1)(2-q)z_0 - q(q-1) \\
&\qquad + (q(q+1)(2-q)z_0 + q(q-1)^2) - (q-1)(2-q)z_0^{1-q}] \\
&= z_0^{q-1}[(q-1)(q+1)(2-q)z_0 \\
&\qquad + q(q-1)(q-2) - (q-1)(2-q)z_0^{1-q}] \\
&= (q-1)(2-q)z_0^{q-1}[(q+1)z_0 - q - z_0^{1-q}] \\
&= (q-1)(2-q)z_0^{q-1} f_1(z_0) \\
&\leq 0 \quad \text{by (A.2) and that } 1 < q < 2.
\end{aligned}
$$

Thus,

(A.3) $$\sup_{0 \leq z \leq 1} f_2(z) \leq \max[f_2(0), f_2(1), 0] = 0.$$

Now define

$$
\begin{aligned}
f_3(z) &= -(q+1)(2-q)z^2 - 2q(q-1)z + q(q-1) \\
&\quad + 2qz^{3-q} - 2(q-1)z^{2-q}, \qquad z \geq 0.
\end{aligned}
$$

Then

$$
\begin{aligned}
f_3'(z) &= 2z^{1-q}[-(q+1)(2-q)z^q - q(q-1)z^{q-1} \\
&\qquad + q(3-q)z - (2-q)(q-1)] \\
&= 2z^{1-q} f_2(z).
\end{aligned}
$$

By (A.3),

$$f_3'(z) \leq 0 \text{ on } (0,1).$$

Thus,

(A.4) $$\inf_{0 \leq z \leq 1} f_3(z) = f_3(1) = 0.$$

Finally, define

$$
\begin{aligned}
f_4(z) &= -(2-q)z^{q+1} - 2(q-1)z^q \\
&\quad + qz^{q-1} + qz^2 - 2(q-1)z - (2-q), \qquad z \geq 0.
\end{aligned}
$$

Then

$$f_4'(z) = z^{q-2} f_3(z) \geq 0 \quad \text{for } z \in (0,1), \text{ by (A.4).}$$

This implies

(A.5) $$\sup_{0 \leq z \leq 1} f_4(z) \leq f_4(1) = 0.$$

Routine computation shows

$$f'(z) = \frac{f_4(z)}{z^{3-q}(1-z^{q-1})^3} \leq 0 \qquad \text{for } z \in (0,1), \text{ by (A.5)}.$$

Thus, $f$ is decreasing on $(0,1)$, as desired. $\square$

**Acknowledgment.** I thank the referee for a speedy and thorough review of this article. The numerous detailed suggestions for clarification and simplification have greatly improved the exposition.

DEPARTMENT OF MATHEMATICS
TEXAS A&M UNIVERSITY
COLLEGE STATION, TEXAS 77843-3368
USA
E-MAIL: deblass@math.tamu.edu